\documentclass[oneside,english,reqno]{amsart}
\usepackage[T1]{fontenc}
\usepackage[latin9]{inputenc}
\usepackage{units}
\usepackage{amsthm}
\usepackage{amstext}
\usepackage{amssymb}
\usepackage[pdfborder={0 0 0},pdfborderstyle={},colorlinks=true,linkcolor=black,citecolor=black,urlcolor=blue]{hyperref}

\makeatletter
\theoremstyle{plain}
\newtheorem*{thm*}{Theorem}
\newtheorem{thm}{Theorem}
\newtheorem{lem}{Lemma}

\theoremstyle{definition}
\newtheorem*{rem*}{Remark}
\newtheorem*{xca*}{Exercise}

\makeatletter
  \let\my@@font@warning\@font@warning
  \let\@font@warning\@font@info
\makeatother
\usepackage{fourier}
\makeatletter
  \let\@font@warning\my@@font@warning
\makeatother

\makeatother

\newcommand{\prob}{\mathbf{Pr}}
\newcommand{\CC}{{\mathbb C}}
\newcommand{\PP}{{\mathbb P}}
\newcommand{\QQ}{{\mathbb Q}}
\newcommand{\ZZ}{{\mathbb Z}}
\newcommand{\Cc}{{\mathcal C}}

\usepackage[english]{babel}
\title{On common roots of random Bernoulli polynomials}
\author{Gady Kozma
\and Ofer Zeitouni}
\begin{document}
\begin{abstract} We prove that with high probability,
	$d+1$ random Bernoulli polynomials in $d$ variables of
degree $n$ ($n\to\infty$) 
do not possess a common root. 
\end{abstract}
\maketitle
\section{Introduction}
We consider in this paper systems of random polynomials in $d$ variables
with independent
Bernoulli coefficients, and study whether they possess a common root.
Specifically, with $\vec{j}_d=(j_1,\ldots,j_d)$, $j_i\in \mathbb{Z}_+$ and
$|\vec{j}_d|=\sum j_i$, let
$\{\epsilon_{ \vec{j}_d}\}$ be a family of i.i.d. Bernoulli, 
$\pm1$-valued random variables.
Set
$x^{ \vec{j}_d}=\prod_{i=1}^d x_i^{j_i}$.
We call the following polynomial
$$ P(x_1,\ldots,x_d)=\sum_{ \vec{j}_d: |\vec{j}_d|\leq n} \epsilon_{
\vec{j}_d} x^{ \vec{j}_d} $$
a {\it random polynomial in \
$d$ variables and degree $n$}. Our main goal
in this paper is to prove the following.
\begin{thm}
	\label{theo-1}
Let $P_{1},\dotsc,P_{d+1}$ be $d+1$ independent
random 
polynomials in $d$ variables and degree \
$n$. Let
$$ p(n,d)=\prob(\exists x \in \mathbb{C}^d: 
P_i(x)=0, i=1,\ldots, d+1)$$
denote the probability that the $P_i$
have a common zero.
Then, there exists a constant $c(d)<\infty$ such that, for all $n$ positive
integer,
\begin{equation}
	\label{eq-1}
p(n,d)\leq  c(d)/n\,.
\end{equation}
\end{thm}
In particular, with probability approaching $1$ as $n\to \infty$, there does
not exist a common zero for the polynomials $P_i$, $i=1,\ldots,d+1$, an
intuitively obvious but otherwise non trivial fact. We remark that the
result would be trivial if the distribution of the coefficients of the
polynomials were to have a continuous distribution --- in this case
the probability would simply be 0, for all $d$ and all $n$. The point
about the result is the discreteness of the coefficients, and we chose
Bernoulli as the simplest example.

Another simple point to note is that it is important that the
distribution has no atom at 0. Indeed, if we were to take Bernoulli
variables taking values 0 and 1 (rather than $\pm 1$), then there
would be probability $2^{-d-1}$ that $(0,\dotsc,0)$ is a common root:
all you need for this event is that the constant coefficient of all
$d+1$ polynomials would be 0, which happens with probability
$2^{-d-1}$ independently of $n$.

The structure of the paper is as follows. In the next section we consider
the case $d=1$. In section \ref{sec-break} we break the event of existence of a common zero according to the type of zero, i.e. according to whether the common zero has at least one zero component, whether it
satisfies a relation determined by 
two monomials (a ``dunomial''), or whether it satisfies neither condition; we
handle the first case by a dimension reduction argument and the last
by a projectivization argument. The dunomial case (and completion of the proof
of Theorem \ref{theo-1}) are presented in section \ref{sec-dunom}, where
a version of Hal\'{a}sz' theorem plays a decisive role.

{\it Convention}. Throughout, $C$ denotes a constant independent of $d$ and $n$ that may change from line to line; 
$c=c(d)$ denote constants
that depend only on $d$ but may change from line to line. On the other hand, 
constants of the form $c_i(d)$  depend on $d$ only
and do not change from line to line.

\section{The one-dimensional case}
In this section we will treat the one-dimensional case. 
Recall that all our random polynomials have $\pm 1$ Bernoulli coefficients.
We will prove

\begin{thm}\label{thm:d=1}
Let $P_1,P_2$ be two independent random polynomials in $1$ variable 
of degree
$n$. Let
\[
p(n)=\prob(\exists x\in\CC : P_1(x)=P_2(x)=0).
\]
If $n$ is even then
\[
p(n)=\Big(\frac{4}{\pi}+o(1)\Big)\frac{1}{n}
\]
and if $n$ is odd then $p(n)\leq c(d) n^{-3/2}$.
\end{thm}
We note that 
the techniques of the next sections apply to the one-dimensional case
unchanged, and could yield an upper bound 
in theorem \ref{thm:d=1} of the form $c(d) n^{-1/2}$. However,
in the one dimensional case there 
are a few additional tools that  
yield a more precise result. 
As we will see, the techniques allow, in principle, to get an
asymptotic series, but we will not go that far in the direction of
extra-precise results.

\begin{proof} The proof relies on two observations.
\begin{itemize}
\item Any solution $\xi$ of $P_1$ must satisfy $\frac{1}{2} < |\xi| < 2$.
\item Any solution $\xi$ of $P_1$ must be an algebraic integer.
\end{itemize}
The first observation is obvious: if $|\xi|\ge 2$ then the highest term
$\epsilon_n \xi^n$ dominates all the others and the sum cannot be zero, if
$|\xi|\le\frac{1}{2}$ then the lowest term $\epsilon_0$ dominates all
others. The second observation is by definition: an algebraic integer
is defined as a number satisfying a monic polynomial, i.e.\ a polynomial with integer
coefficients and highest coefficient equal to $1$. Nevertheless, there
is an algebraic fact in the background of the definition, which we
will use

\medskip
\noindent {\em A number is a root of a monic polynomial if and only if it is
  the root of an irreducible monic polynomial.}

\medskip
\noindent See any standard textbook on algebraic number theory, e.g.\ \cite[Page 93]{AW}. These two observations allow us to
classify all potential solutions of low algebraic order. For example,
if $\xi$ is a rational solution then its irreducible polynomial is
$x-\xi$ and since it is monic, $\xi$ must be an integer, and by the
first observation it must be $\pm 1$. If the irreducible polynomial of
$\xi$ is of degree $2$, it must be $x^2 - ax -b$. But because it is
irreducible, the other solution $\xi'$ is also a solution of $P_1$, so
it must also be between $\frac{1}{2}$ and 2. But $a=\xi+\xi'$ and
$b=\xi\xi'$ so both are between $4$ and $-4$. We get that we only need
to examine a finite collection of numbers (naively 98, since $a$ has 7
possibilities, $b$ has 7 possibilities, and each polynomial has two
roots --- this number can be reduced easily, but this is not important
at this step). The same argument gives
\begin{lem}\label{lem:finite}
For every $\ell$ there are only finitely many numbers
whose irreducible polynomial has degree $\ell$ that can be roots of a
polynomial (of arbitrary degree) with coefficients $\pm 1$.
\end{lem}
Next let us recall the so-called
``sharp inverse Littlewood-Offord theorem'' of Tao and Vu
\cite[theorem 1.9]{TV10}, which we now quote almost literally:

\medskip
\noindent {\em Let $A$, $\delta>0$ and let $(\xi_1,\dotsc,\xi_n)$ be
  complex numbers such that
\[
\prob\Big(\sum_{i=0}^n \epsilon_i\xi_i=0\Big)\ge n^{-A}.
\]
Then there exist a symmetric generalized arithmetic progression (all
of whose elements are distinct) of rank $B\le 2A$ and volume $\le
C(A,\delta) n^{A-B/2+C(A)\delta}$ which contains all but
$C(A,\delta)n^{1-\delta}$ of the $\xi_i$ (counting multiplicity).
}

\medskip
\noindent The value of $\delta$ will play no rule, so we set it to
$\frac{1}{2}$. We apply this theorem with $\xi_i=\xi^i$ and get that if
$\prob(\sum_{i=0}^n
\epsilon_i\xi^i)\ge n^{-A}$ then most $\xi_i$s  must be
contained in a generalized arithmetic progression of rank $B\le 2A$.
More precisely,
there exists some $\gamma_1,\dotsc,\gamma_B\in\CC$ such that
we have $\xi^i=\sum_j n_{i,j}\gamma_j$ with $n_{i,j}$
integers for at least $n- 
Cn^{1-\delta}=n-C\sqrt{n}$ indices $i$. Therefore, for $n$ sufficiently
large,
$n>n_0(A,\delta)$, it must hold  for some $B+1$ consecutive $i$s, call them
$i,\dotsc,i+B$. But these $B+1$ vectors of coefficients $n_{i,j}$ must be
dependent over the rationals $\QQ$, so $\xi^i,\dotsc,\xi^{i+B}$ must
be dependent over $\QQ$, which means that
$\xi$ satisfies a polynomial with rational coefficients of degree $\le
B$. In other words we proved
\begin{lem}\label{lem:TV} For all $A>0$ there exists $n_0(A)$ such that if
	$n>n_0(A)$ and if $\prob(\sum_{i=0}^n 
	\epsilon_i \xi^i=0)\ge n^{-A}$ then $\xi$
  must be of algebraic degree $\le 2A$.
\end{lem}

Let us finish the proof of theorem \ref{thm:d=1}. Using the remarks
before lemma \ref{lem:finite} we write
\begin{align}
\lefteqn{\prob(\exists x : P_1(x)=P_2(x)=0) =}\qquad\qquad\nonumber\\
 \;&\prob(P_1(1)=P_2(1)=0) \nonumber\\
& + \;\prob(P_1(-1)=P_2(-1)=0) \nonumber\\
& - \;\prob(P_1(1)=P_1(-1)=P_2(1)=P_2(-1)=0) \nonumber\\
& + \;O(\prob(\exists x \textrm{ of algebraic degree }\in \{2,3,4,5\} :
P_1(x)=P_2(x)=0)\nonumber\\
& \qquad+
\prob(\exists x \textrm{ of algebraic degree }> 5 : P_1(x)=P_2(x)=0))
\nonumber\\ 
= &\; I + II + III + O(IV + V)\label{eq:I-V}
\end{align}

The estimate of the first three terms is straightforward. The first
term ($I$) is exactly the probability that a random walk on $\ZZ$
returns to 0 at time $n$, squared, since we need both $P_1$ and $P_2$
to be zero. This can be estimated by Stirling's formula and we get
\[
I=II=\begin{cases}\left(\frac{2}{\pi}+o(1)\right)\frac{1}{n} & n\textrm{
  is even}\\
0 & n\textrm{ is odd.}
\end{cases}
\]
The third term $III$ is the probability that $P_1$ and $P_2$ 
both vanish at $\pm 1$.
This probability vanishes
when $n$ is odd and, when $n$ is even,
it equals the probability that both $\sum_{i=0}^{n/2}
\epsilon_{2i}=0$ and 
$\sum_{i=0}^{n/2}
\epsilon_{2i+1}=0$; those events are independent, so $|III|\le Cn^{-2}$
($III$ is negative).

For the term $IV$ we use lemma \ref{lem:TV} with $A=\nicefrac{3}{4}$ and
see that any $\xi$ such that $\prob(\sum\epsilon_i\xi^i=0)>n^{-3/4}$ must be
rational, so does not contribute to $IV$. So we get that any $\xi$
with algebraic degree $\in\{2,3,4,5\}$, 
\[
\prob(P_1(\xi)=P_2(\xi)=0)\le \left(n^{-3/4}\right)^2=n^{-3/2}.
\] 
By lemma \ref{lem:finite} there are only finitely many $\xi$ which we
need to consider, so $IV\le Cn^{-3/2}$.

Finally, for the term $V$ we fix $P_1$. It has (at most) $n$ different
roots. For each one we ask what is the probability that it is also a
root of $P_2$? We use lemma \ref{lem:TV} with $A=\nicefrac{5}{2}$ and
get that, for $n$ sufficiently large, any $\xi$ such that $\prob(P_2(\xi)=0)\ge n^{-5/2}$ has
algebraic degree $\le 5$ so does not contribute to $V$. So we can
write
\begin{align*}
V&=\prob(\exists\textrm{a root of $P_1$ with algebraic degree $>5$
  which is also a root of $P_2$}) \\
&\le n\max\{\prob(P_2(\xi)=0):\xi
\textrm{ with algebraic degree > 5}) \le n \cdot n^{-5/2} = n^{-3/2}.
\end{align*}
Plugging the estimates for $I$--$V$ into (\ref{eq:I-V}) finishes the
proof of theorem \ref{thm:d=1}.
\end{proof}

\section{Breakup into cases}
\label{sec-break}
The proof of Theorem \ref{theo-1} goes by induction on the number of variables $d$. Before starting
it, we introduce some notation and
prove auxiliary lemmas.

In the sequel, for a collection of polynomials $Q_1,\ldots,Q_\ell$, we write
$Z(Q_1,\ldots,Q_\ell)$ for their common zeros, i.e. the algebraic set
determined by this collection. This algebraic set
may be reducible. We denote
\begin{description}
\item[$Z_\infty(Q_1,\dotsc,Q_\ell)$] The union of all irreducible components of $Z$ with
  dimension $>0$.
\end{description}
(we call this ``$Z_\infty$'' because it contains infinitely many
points, if non-empty). For the definition of irreducible components of
an algebraic variety, see any standard textbook,
e.g.~\cite{S74}. 

Next divide $Z=Z_1\cup Z_2\cup Z_3$ (these are not directly related to
$Z_\infty$ --- we hope the reader will not
be too confused by the somewhat inconsistent use of the subscript), as follows:
\begin{description}
\item[$Z_1(Q_1,\dotsc,Q_\ell)$] The elements $(x_1,\dotsc,x_d)\in Z$ which satisfy a monomial, or in other
  words, that (at least) one of the $x_i$ is zero.
\item[$Z_2(Q_1,\dotsc,Q_\ell)$] The elements $(x_1,\dotsc,x_d)\in Z$ which do not satisfy a monomial but
  do satisfy a
  dunomial\footnote{Binomial might have been a better term, but is
    already taken in the literature} of degree at most $n$,
    i.e. such that for some
  $\vec{\alpha}\ne\vec{\beta}$ with $\sum \alpha_i\leq n$, 
  $\sum \beta_i\leq n$,
\[
\prod_{i=1}^d x_i^{\alpha_i}\pm\prod_{i=1}^d x_i^{\beta_i}=0
\]
\item[$Z_3$] The elements of $Z$ which satisfy neither a monomial nor
  a dunomial.
\end{description}

Applying this to $\ell$ random polynomials $P_1,\dotsc,P_\ell$ of
degree $n$ in $d$
variables we define the following
corresponding probabilities
\[
p_i(n,d,\ell)=\prob(Z_i(P_1,\dotsc,P_\ell)\ne\emptyset)\qquad i\in\{1,2,3,\infty\}
\]
We first estimate $p_\infty$ --- we believe this is the most
interesting estimate in the proof (it definitely took us longest to discover).
\begin{lem}
	\label{lem-1}
	For any $d\geq 2$
	and  all $n$
	positive integer,
	\begin{equation}
		\label{eq-2}
		p_\infty(n,d,\ell)\leq d  p(n,d-1,\ell)\,.
	\end{equation}
\end{lem}
\begin{proof}
	Let $\Cc$ be an arbitrary irreducible component
	(of dimension necessarily $\geq 1$) of 
	$Z_\infty(P_1,\ldots,P_\ell)$.
	 We examine
$\Cc$ in the $d$ dimensional projective space $\mathbb{P}^{d}$ i.e.\ add
a $d+1^{\textrm{st}}$ variable and homogenize by multiplying each
monomial $x_{1}^{\alpha_{1}}\dotsb x_{d}^{\alpha_{d}}$ by 
$x_{d+1}^{n-\sum\alpha_{i}}$
so we get a system of homogeneous polynomials of degree $n$ in $d+1$
variables. 
A nice feature is that there is no difference between the
added variable and the old ones --- our polynomials are \[
\sum_{\vec{\alpha}_{d+1}:
\sum\alpha_{i}=n}\epsilon_{\vec{\alpha}_{d+1}}
\prod_{i=1}^{d+1}x_{i}^{\alpha_{i}}.\]
Clearly the set of zeros of the homogenized system has a component
of dimension $\ge2$
because \[
\{(\lambda x_{1},\dotsc,\lambda x_{d},\lambda):(x_{1}\dotsc,x_{d})\in \Cc,
\lambda\in\mathbb{C}\}\]
are all zeros. Hence the dimension of $\Cc$ as a projective variety
(denoted now as $\widetilde{\Cc}$)
is $\ge1$.

We now apply the projective dimension theorem \cite[\S 1, theorem 7.2]{H77}
to see that $\widetilde{\Cc}$ intersects with the plane $x_{1}=0$. 
Call this intersection
$(\mu_{1},\dotsc,\mu_{d+1})$. By definition they cannot be all zero
--- this is not a legal point in the projective space. So let $k>1$
satisfy that $\mu_{k}\ne0$, and in this case we may assume $\mu_{k}=1$
and remove $\mu_{k}$ from our equations. We are left with $d-1$
variables: $\mathcal{V}=\{2,\dotsc,d+1\}\setminus\{k\}$. So we get
that a system of $\ell$ independent polynomials in $d-1$ variables
\[
\sum_{\vec{\alpha}_{d-1}:\sum_{i\in\mathcal{V}}\alpha_{i}\le n}\epsilon_{\vec{\alpha}_{d-1}}\prod_{i\in\mathcal{V}}x_{i}^{\alpha_{i}}\]
has a common zero. In other words, if we denote the event that this
system of equations has a common zero by $E_k$, then the conclusion is
the event that
$Z_\infty(P_1,\dotsc,P_\ell)\ne\emptyset$ implies $E_1\cup\dotsb\cup E_d$.
By definition each of the $E_k$ has probability $p(n,d-1,\ell)$.
We need to count over $k$, which has $d$ possibilities, so we get
that $p_\infty(n,d,\ell)\leq dp(n,d-1,\ell)$.
\end{proof}


We proceed with estimates of $p_i$ for finite $i$. The obvious one is
\begin{lem}\label{lem:p1}
\[
p_1(n,d,\ell)\le d p(n,d-1,\ell).
\]
\end{lem}
\begin{proof} If $x_i=0$ for some $i$ we can throw it and all terms
  containing it and we get exactly $p(n,d-1,\ell)$. The term $d$ comes from counting over the $i$.
\end{proof}
The second easiest one is
 \begin{lem}
	 \label{lem-2}
	 For all $d\ge 1$ there exists a constant $c_1=c_1(d)$ independent of
	 $n$ such that for all $n$ positive
	 integer and all $\ell$,
	 \begin{equation*}
		  p_3(n,d,\ell)\leq p_\infty(n,d,\ell-1)+c_1 n^{-d/2}\,.
	  \end{equation*}
  \end{lem}
  \begin{proof}
Examine $P_1,\dotsc,P_{\ell-1}$. The event that
$Z(P_1,\dotsc,P_{\ell-1})$ has a component of dimension $\ge 1$ we
push into the $p_\infty$ term, so we may assume that all components
are points. By Bezout's theorem \cite{H83, T}, the cardinality of
$Z(P_1,\ldots,P_{\ell-1})$ is at most $n^d$. This of course applies
also to $Z_3$ which is a subset of $Z$.

We now add the last polynomial $P_\ell$. We need to ask, for every
$\vec{x}\in Z_3(P_1,\dotsc,P_{\ell-1})$, what is the
probability that $P_\ell(\vec{x})=0$? We apply the
S\'ark\H{o}zy-Szemer\'edi theorem \cite{SS65} which states that for any fixed $\xi\in
\CC^{m}$ with all $\xi_i$ different,
\[
\prob\bigg(\sum_i\epsilon_i\xi_i=0\bigg)\le
cm^{-3/2}.
\]
For $\vec{x}\in Z_3$ we know that they
satisfy no dunomial, hence the vector 
\[
\xi_{\vec{j}}=\prod_{i=1}^d x_i^{j_i}
\]
(which lives in $\CC^m$ for $m$ being the number of possible choices
of $\vec{j}$, so $m\approx n^d$) has all entries distinct. Hence 
\[
\prob\Big(\sum_{\vec{j}}\epsilon_{\vec{j}}\prod_{i=1}^d
x_i^{j_i}=0\Big)\le cm^{-3/2}
\]
Since $m \ge c(d)n^d$ we get that for each $\vec{x}\in
Z_3$ we have $\prob(P_\ell(\vec{x})=0)\le c_1 n^{-3d/2}$. Summing over
all $\vec{x}$ and using the information gathered from Bezout's theorem
finishes the lemma. 
 \end{proof}

\begin{rem*} Using lemma \ref{lem-1} and the idea of lemma
  \ref{lem-2} (with
  the S\'ark\H{o}zy-Szemer\'edi theorem replaced by Erd\H{o}s' theorem
  \cite{E45}), one can show that
  a system of $3d-1$ random polynomials in $d$ variables of degree $n$
  does not have
  a common root, with high probability. 
(Recall that Erd\H{o}s' theorem states that if all $\xi_i$ are
non-zero, then $\prob(\sum_{i=0}^n\epsilon_i\xi_i=0)\le Cn^{-1/2}$.)
The other arguments in the paper, and most notably the use of dunomials
and the Hal\'{a}sz theorem, are needed in order to reduce the number of required
polynomials from $3d-1$ to $d+1$.
\end{rem*}

\section{Dunomial analysis and proof of Theorem \ref{theo-1}}
\label{sec-dunom}
Lemma \ref{lem-2} gives a handle on analyzing the set $Z(P_1,\ldots,P_{d+1})$,
away from
the set of points that are zeros of a dunomial, and do not possess zero 
coordinates. To be able to carry an induction step and provide a proof of
Theorem \ref{theo-1}, we thus need to consider such points.

Let
$$ D(x)=\prod_{i=1}^d x_i^{\alpha_i} \pm 
\prod_{i=1}^d x_i^{\beta_i}\, $$
be a dunomial of $d$ variables and degree less than or equal to
$n$. Define the {\em order} of $D$ to be 
\[
|D|=\sum_{i=1}^d |\alpha_i-\beta_i|
\]

For $\vec{x}\in (\mathbb{C}\setminus \{0\})^d$, let 
\[
r(\vec{x})=\min\{|D|:D(\vec{x})=0\}
\]
i.e.~the minimal order among all dunomials satisfying $\vec{x}$, or
$\infty$ if none satisfy it.
The following lemma is simple but crucial.
\begin{lem}
	\label{lem-3}
	There exists a constant $c_2(d)$ so that for any $n$ positive integer
	and  $x\in (\mathbb{C}\setminus \{0\})^d$, the number $R_n(x)$ of
        dunomials $D$ of degree $\le n$ satisfied by $x$ has
\begin{equation*}
	R_n(x) \leq 
	c_2(d)\frac{n^{2d}}{r(x)^{d}}\,.
\end{equation*}
\end{lem}
\begin{proof}
        Fix some $\vec{\alpha}$, and assume $\vec{\beta}$ satisfies
        that $\prod x_i^{\alpha_i}=\pm \prod x_i^{\beta_i}$. If we
        have for some $\vec{\gamma}$ that also $\prod
        x_i^{\alpha_i}=\pm\prod x_i^{\gamma_i}$ then by definition we
        must have $|\vec{\beta}-\vec{\gamma}| \ge r(x)$. Thus, $R_n(x)$
	is bounded by the total number of dunomials of degree $n$ (which is
	$\le c(d) n^{2d}$), divided by the minimal number of integer points in a
	ball of $\ell_1$ radius $r(x)$. Since the latter is 
	bounded by a constant (depending on $d$) multiple of
	$r(x)^{d}$, the lemma follows.
\end{proof}

\begin{lem}
	\label{cor-1}
	For any $d\geq 1$ there exists a constant $c_3(d)$ such that
	for all $n$, $\ell$ positive integers,
	\begin{equation}
		\label{eq-7}
		p_2(n,d,\ell)\leq d^3p(n,d-2,\ell-2) + \frac{c_3(d) (\log n)^{(2-d)_+}}{n}
	\end{equation}
(where for $d\le 2$ we use the convention that $p(n,d-2,\ell-2)=0$ in
        the overdetermined case $\ell>d$)\end{lem}
\begin{proof}
	For any dunomial $D$, consider the zero set
	$Z=Z(D,P_1,\ldots,P_{\ell-2})$. Since we are interested in
        $p_2$, we may assume $D$ is reduced i.e.~no $x_i$ appears on both
sides, and we will make this assumption throughout the proof. There are two events to consider
        (corresponding to the two term in (\ref{eq-7})), one where $Z$
        has components of dimension $\ge 1$, and the other where it
        does not.

Let us start with the case that $Z$ does have components of dimension
$\ge 1$. In this case we do not need to know the exact value of the
dunomial $D$ --- we only need to know which $x_i$ appear on the two
terms. We now repeat the analysis of lemma \ref{lem-1} namely embed one
irreducible component $\Cc$ of dimension $\ge 1$ into the projective space $\PP^d$ and use the
projective dimension theorem. This time, however, we do not intersect
$\widetilde{\Cc}$ necessarily with $x_1=0$, but we intersect it with
some $x_i=0$ for some $x_i$ that appears in the dunomial $D$. The
intersection is still non-empty, and of course, if $x_i=0$ and
$D(\vec{x})=0$ then at least one other $x_j$ (appearing in the other
term of $D$ i.e.~in the term not containing $x_i$) must also be
zero (here we use that $D$ is reduced). Now return to the affine setting $\CC^d$ as in lemma
\ref{lem-1}, i.e.~find some $k$ such that the solution $\vec{x}$ has
$x_k\ne 0$ and set $x_k=1$. Recall that in
lemma \ref{lem-1} we defined events $E_k$ that the system one gets by
setting $x_1=0$ and $x_k=1$ has a common root. Here we need instead
events $E_{i,j,k}$ that the system one gets by setting $x_i=x_j=0$ and
$x_k=1$ has a common root. But the conclusion is the same: if for some
reduced $D$ one has that $Z(D,P_1,\dotsc,P_{\ell-2})$ has
a component of dimension $\ge 1$ then necessarily one of the
$E_{i,j,k}$ happened, and each one has probability
$p(n,d-2,\ell-2)$. This explains the first term in (\ref{eq-7}).

Now assume $Z$ is finite. 
Since $D$ is reduced, the order of $D$ is also its degree, and by Bezout's theorem
	$|Z_2(D,P_1,\ldots,P_{\ell-2})|\leq |D|n^{d-1}$. 

Next fix some $\vec{x}\in Z$. We now apply a
        strengthening of the S\'ark\H{o}zy-Szemer\'edi theorem due to
        Hal\'asz \cite{hal}. Let us recall Hal\'asz' theorem. It states that for
        any $\vec{\xi}\in\CC^m$,
\[
\prob\Big(\sum_j\epsilon_j\xi_j=0\Big)\le\frac{R}{m^{-5/2}}
\]
where $R$ is the number of couples $j,k$ (not necessarily different)
such that $\xi_j=\pm\xi_k$. In our case (i.e.~when $\xi_{\vec{j}}=\prod
x_i^{j_i}$), $R$ is exactly given by lemma
\ref{lem-3}: $R$ is $n^{d}+$ the number of dunomials satisfied by $x$
(the $n^d$ term corresponds to the trivial couples $j=k$). Hence $R
\le cn^{2d}/r(x)^d$ as the $n^d$ term is always smaller and hence can be
incorporated into the constant.
 
Thus we get, for any point $\vec{x}\in Z$ with $r(\vec{x})=r$,
	\begin{equation*}
		\prob(P_{\ell-1}(x)=P_{\ell}(x)=0)
		\leq c(d) \left( \frac{n^{2d}}{r^d}\cdot \frac{1}{n^{5d/2}}
		\right)^2\leq c(d)\frac{1}{r^{2d} n^{d}}\,.
	\end{equation*}
	Since there
	are at most $c(d) r^{d-1}$ reduced forms of dunomials of order
	$r$, we  conclude that
	\begin{multline*}
		\prob(\exists D : 
                |Z(D,P_1,\dotsc,P_{\ell-2})|<\infty\textrm{ and }
                  \{\exists \vec{x}\in Z : r(\vec{x})=|D| \textrm{ and }
                P_{\ell-1}(\vec{x})=P_\ell(\vec{x})=0\}) \\
\leq c(d) \sum_{r=1}^{2n}  r^{d-1}\cdot \underbrace{n^{d-1}r}_\textrm{Bezout}\cdot\frac{1}{
		r^{2d}n^d}
                \le c_3(d)\begin{cases}
                  \frac{\log n}{n}& d=1\\
                  \frac{1}{n}& d\ge 2\,.
                \end{cases}
	\end{multline*}
	The conclusion follows, since if
        $Z_2(P_1,\dotsc,P_\ell)\ne\emptyset$, either a $D$ and an $x$
        as above exist, or a $D$ exists such that
        $|Z(D,P_1,\dotsc,P_{\ell-2})|=\infty$.
\end{proof}

We can now provide the following.
\begin{proof}[Proof of theorem \ref{theo-1}] The proof is by
  induction on $d$. The case $d=1$ is done by theorem \ref{thm:d=1}. Recall that $p(n,d,\ell)$ is the probability that a
  system of $\ell$ random polynomials in $d$ variables of degree $n$
  has a common root. We write
\begin{align*}
p(n,d,d+1)&=p_1(n,d,d+1)+p_2(n,d,d+1)+p_3(n,d,d+1)\le\\
\intertext{Using lemma \ref{lem:p1} to estimate $p_1$, lemma
  \ref{cor-1} to estimate $p_2$ and lemma \ref{lem-2} to estimate $p_3$,} 
&\le dp(n,d-1,d+1) + d^3p(n,d-2,d-1) + \frac{c}{n} + p_\infty(n,d,d)+c
n^{-d/2}\le\\
\intertext{Using lemma \ref{lem-1} to estimate $p_\infty$,}
&\le dp(n,d-1,d+1) + d^3p(n,d-2,d-1) + \frac{c}{n} + d p(n,d-1,d)\le\\
\intertext{and inductively}
&\le \frac{c(d)}{n}.
\end{align*}
(it is also possible to avoid using theorem \ref{thm:d=1}, and
estimating the 1 dimensional case using these tools. This
will give $p_2(n,1,2)\le (c\log n)/n$ and $p_3(n,1,2)\leq c/\sqrt{n}$, so the
overall result will be that $p(n,1,2)\le c/\sqrt{n}$ and the same
estimate will pass inductively to all $p(n,d,d+1)$.)
\end{proof}

\subsection*{Acknowledgements}
Lots of people have helped us with the algebraic parts of the
proof. Special thanks go to Uri Bader, Amos Nevo,
Dmitry Gourevich, Avraham Aizenbud, Steve Gelbart
and Sergei Yakovenko.
We thank
Ron Peled and Zeev Rudnik for interesting discussions (unrelated to
the algebraic part). Both authors were supported by their respective Israel
Science Foundation grants.

\end{document}